\documentclass[leqno,11pt]{amsart}

\usepackage{amsmath,amsthm}
\usepackage{amsfonts,amssymb}
\newtheorem{theorem}{Theorem}
\newtheorem{lemma}[theorem]{Lemma}
\newtheorem{proposition}[theorem]{Proposition}
\newtheorem{corollary}[theorem]{Corollary}

\def \sgn{{\rm sgn }}
\def \pr{{\mathbb P}}
\def \Pm{{\bf P}}
\def \meanm{{\bf E}}

\def \Par{W}
\def \sgn{{\rm sgn }}
\newcommand{\bb}{\mathbb}

\def \Z {{\bb Z}}
\def \S {{\mathcal S}}
\def \id {{\rm id}}
\def \R {{\bb R}}
\def \invpi{\widetilde{\pi}}

\begin{document}

\title[Series Jackson networks and non-crossing probabilities]{Series Jackson networks \\and non-crossing probabilities}
\author{A. B. Dieker}
\address{Georgia Institute of Technology,
Stewart School of Industrial and Systems Engineering,
765 Ferst Drive NW, Atlanta GA 30332-0205, USA}
\email{ton.dieker@isye.gatech.edu}
\author{J. Warren}
\address{University of Warwick, Department of Statistics, Coventry, CV4 7AL, United Kingdom}
\email{j.warren@warwick.ac.uk}
\date{\today}
\begin{abstract}
This paper studies the queue length process in series Jackson networks
with external input to the first station only.
We show that its Markov transition probabilities can be written as a finite sum of non-crossing probabilities,
so that questions on time-dependent queueing behavior are translated to questions on
non-crossing probabilities. 
This makes previous work on non-crossing probabilities
relevant to queueing systems and allows new queueing results to be established.
To illustrate the latter, we prove that the relaxation time (i.e., the reciprocal of the `spectral gap')
of a positive recurrent system equals the relaxation time of
an M/M/1 queue with the same arrival and service rates as the network's bottleneck station.
This resolves a conjecture of Blanc~\cite{blanc:relaxation1985}, which he proved for two queues in series.
\end{abstract}

\maketitle

\section{Introduction}
The queue length process $\{Q(t):t\ge 0\}$ of an M/M/1 queue with $Q(0)=0$ admits the
sample-path representation
\[
Q(t)= \sup_{0\le s\le t} \left[ N_0(t)-N_1(t)-(N_0(s)-N_1(s))\right],
\]
where both $N_0$ and $N_1$ are Poisson counting processes, $N_0$ being the arrival process and $N_1$
the process counting the number of departures and unused services.
Using time-reversal we find that $Q(t)$ given $Q(0)=0$ equals $\sup_{0\le s\le t} [N_0(s)-N_1(s)]$
in distribution, and in particular that $\pr(Q(t)=0|Q(0)=0)$ can be expressed as a non-crossing
probability
\begin{equation}
\label{eq:onedimcase}
\pr(Q(t)=0|Q(0)=0)=\pr(T>t),
\end{equation}
where $T = \inf\{t\ge 0: N_1(t) < N_0(t)\}$.

This paper generalizes (\ref{eq:onedimcase}) to a network setting,
partially relying on the combinatorial techniques and symmetric functions
discussed in \cite{diekerwarren:determinant2008}.
The networks we study are series Jackson networks with external input 
to the first station only.
We express Markov transition probabilities arising from these networks
as a finite (weighted, signed) sum of non-crossing probabilities for multidimensional Poisson processes.

Time-dependent results for queue lengths in Jackson networks are mostly focused on the case of two stations,
e.g., Blanc~\cite{blanc:relaxation1985} and Baccelli and Massey~\cite{baccellimassey:transient1988}.
The spectral techniques from the work of Kroese, Scheinhardt, and
Taylor~\cite{kroesescheinhardttaylor:spectral2004} seem to be relevant in a time-dependent context,
even though they are developed for the stationary distribution.
O'Connell~\cite{oconnell:pathRS2003} obtains results for series Jackson networks
with an arbitrary number of stations, and these results are most closely related to this paper.
Indeed, as in \cite{oconnell:pathRS2003}, the combinatorial mechanism underlying our results
is the Robinson-Schensted-Knuth (RSK) correspondence.

It is well-known that, through the RSK correspondence,
departures from queues in series are related to non-crossing Poisson processes.
For instance, O'Connell and Yor~\cite{oconnellyor:noncolliding2002} prove that the
cumulative departure process of the last queue is equal in distribution to
the smallest component of Poisson processes conditioned not to cross.
In turn, the latter process is closely related to the Charlier random-matrix ensemble.
Our representations for transition probabilities as a sum of non-crossing probabilities
have a different flavor.
In fact, our arguments do not exploit the RSK correspondence explicitly, but we use a recent
result of R\'akos and Sch\"utz~\cite{rakosschutz:bethe2006} as a starting point.
Their formula is a generalization of the Sch\"utz formula for the transition probabilities
in the totally asymmetric simple exclusion process (TASEP).

The connection between queues in series and non-crossing Poisson processes studied here
can possibly be used to obtain new results on queueing networks that are
beyond the scope of the present paper.
For instance, the recent work of Pucha{\l}a and Rolski~\cite{puchalarolski:asymptotics2007}
on non-crossing probabilities suggests that the `sharp' asymptotic behavior (as $t\to\infty$) 
of the queue-length transition probabilities may be found; 
we discuss this in somewhat more detail at the end of Section~\ref{sec:relaxation}.
It may also be possible to find
a `spectral' representation for queueing probabilities analogous to
the M/M/1 spectral (or integral) representation, see, e.g., \cite{abatewhitt:spectral1988},
from a `spectral' representation for the transition probabilities of non-crossing Poisson processes.

\subsection*{Convergence to stationarity}
Representing Markov transition probabilities as a sum of non-crossing probabilities allows us to
examine the speed at which the multi-dimensional queue length vector $Q(t)$
converges to its steady-state vector $Q(\infty)$.
More precisely, given that all queues are initially empty,
we investigate the probability that the system is empty at time $t$ for large $t$.
Using a large-deviation result for non-crossing probabilities,
we determine the relaxation time for the queueing network,
i.e., the `asymptotic' time required for this probability to decrease to $1/e$.
For many Markovian transition probabilities (including reversible kernels on a finite state space), the
relaxation time is the reciprocal of the spectral gap.
We show that the network's relaxation time equals the relaxation time of
an M/M/1 queue with the same arrival and service rates as the bottleneck station, see (\ref{eq:relaxationtime}).
This resolves a conjecture of Blanc~\cite{blanc:relaxation1985}, which he proved for two queues in series.
For work on relaxation times for non-Markovian processes, we refer to the recent work of Glynn,
Mandjes, and Norros~\cite{glynnmandjesnorros:convergence2008}.

\subsection*{The model}
We consider a queueing system with several single-server queues in series,
each with an unlimited waiting space.
Customers arrive according to a Poisson process at the first queue.
They join station $i+1$ after receiving service at station $i$, and leave the
system after being served at the last station.
The service discipline in each of the stations is first-in-first-out.
All service times are exponentially distributed with a parameter only depending on the station,
as well as mutually independent and independent of the arrival process.
It is convenient to imagine that the arrival stream arises
as the output process of an auxiliary zeroth station
with an infinite number of customers;
the service requirements at this station then become the interarrival times.
We suppose that the total number of stations (excluding the auxiliary
station) is $N$, and we let $\nu_i$ be the parameter of the exponential
service distribution at station $i$. In particular, $\nu_0$ is the
arrival rate of customers at the first station.

Two natural Markov processes associated to series Jackson networks are the cumulative number
of departures vector $D$ and the queue length vector $Q$.
For $i=0,\ldots,N$ and $t\ge 0$, let $D_i(t)$ be the total number of departures from the $i$-th
queue in the time interval $[0,t]$. For $i=1,\ldots,N$, write $Q_i(t)$ for the queue
length at station $i$ at time $t$.
The vector $Q$ is readily described in terms of $D$: $Q(t)=(Q_1(0)+D_0(t)-D_1(t),Q_2(0)+D_1(t)-D_2(t),\ldots,Q_{N}(0)+D_{N-1}(t)-D_{N}(t))$.

\subsection*{Organization and notation}
This paper is organized as follows.
Section~\ref{sec:departures} presents results on the cumulative departure process which are principal tools
in our investigations. Section~\ref{sec:identities} derives two identities for the probability of an empty system
at time $t$ given an initially empty system, the second identity being only applicable in the positive recurrent case.
Relaxation times are studied in Section~\ref{sec:relaxation}. In Section~\ref{sec:discussion}, we
discuss how the identities of Section~\ref{sec:identities} change when the underlying assumptions are relaxed.
We review the concepts of combinatorics and symmetric functions that are most relevant to us in Appendix~A.

Throughout, we let the family $\{w_n:n\in\Z\}$ of functions on $\R$ be given by
$w_n(t)=t^n/n! 1_{\{t\ge 0\}}$ for
$n>0$, while $w_0(t)=1_{\{t\ge 0\}}$ and $w_n=0$ for $n<0$. Here and throughout, $1_A$ denotes the indicator
function of the set $A$.
We also define
\[
\Par^N=\{x\in\Z^{N+1}: x_0\ge \ldots \ge x_N\}.
\]
Unless otherwise indicated,
all determinants in the present paper are $(N+1)\times (N+1)$, with the
indices $i$ and $j$ ranging from $0$ to $N$.

\section{The cumulative number of departures}
\label{sec:departures}
This section focuses on the probabilistic behavior of
the cumulative number of departures from each of the queues,
by presenting the transition kernel of this process.
We also express this kernel in terms of the kernel arising
from non-colliding Poisson processes.

A full description of the process $D$ on a probabilistic level is given by the transition probabilities
$\phi^\nu_t$, defined for $d,d'\in\Par^N$ and $t>0$ as
\[
\phi^\nu_t(d,d')=\pr(D(t)=d'|D(0)=d).
\]
R\'akos and Sch\"utz~\cite{rakosschutz:bethe2006}
derive an explicit formula for $\phi^\nu_t$, which we describe next.
They study an interacting particle process known as totally asymmetric simple exclusion process (TASEP),
which is equivalent to the departure process from queues in series.
Both the vector-valued process $D$ and the TASEP particle system are Markov jump processes.
The departure process at the $i$-th queue $D_i$ makes a jump of $+1$ at rate $\nu_i$, corresponding to a service (or arrival if $i=0$), unless $D_i(t)=D_{i-1}(t)$ in which case it cannot
jump because the $i$-th queue is empty. Now consider the `shifted' vector
$\tilde D$ where $\tilde D_i(t)=D_i(t)-i$, and note that this process has the same dynamics
as the ordered positions of particles in an $(N+1)$-particle TASEP exclusion process 
with particle $i$ having jump rate $\nu_i$ for $i=0,\ldots,N$.

The next proposition presents the result by R\'akos and Sch\"utz in a form which
closely follows \cite{diekerwarren:determinant2008}.
R\'akos and Sch\"utz~\cite{rakosschutz:bethe2006} prove this proposition by
verifying the Kolmogorov forward equations for $D$,
and an alternative proof can be given using the results from
\cite{diekerwarren:determinant2008}.
For $n\in\Z,t\in\R$ and $0\le i,j\le N$, we write
\begin{equation}
w^{(ij)}_n(t)=
\begin{cases}
\sum_{k=0}^{i-j} (-1)^k e_k^{(ji)}(\nu) w_{n+k}(t)
&\text{ if } j\le i, \\
\sum_{k=0}^\infty h_k^{(ij)}(\nu) w_{n+k}(t)
&\text{ if } i\le j,
\end{cases}
\end{equation}
where $e_k^{(ji)}(\nu)$ and $h_k^{(ij)}(\nu)$ are coefficients which are elementary functions of the vector $\nu=(\nu_0,\ldots,\nu_N)$;
see Appendix~A for explicit formulas. 
In the remainder, for notational convenience, we suppress the arguments $\nu$ of 
the coefficients $e_k^{(ij)}(\nu)$ and $h_k^{(ij)}(\nu)$.
As an aside, we note that $w^{(ij)}_n(t)$ is a shifted multiple orthogonal polynomial of Charlier~\cite{arvesu:multipleorthogonal2003} for $j\le i$.

\begin{proposition}[R\'akos-Sch\"utz~\cite{rakosschutz:bethe2006}]
For $d,d'\in\Par^{N}$, we have for any $t> 0$,
\[
\phi^\nu_t(d,d')=
\prod_{k=0}^N \left[e^{-\nu_k t} \nu_k^{d_k'-d_k} \right]
\det \left\{w^{(ij)}_{d_i'-d_j-i+j}(t)\right\}.
\]
\end{proposition}

It is our next aim to derive a representation for $\phi^\nu_t$ in terms of the transition
kernel of Poisson processes killed on their first crossing time.
To describe these processes in more detail,
consider a measurable space on which a stochastic process
$\{X(t)=(X_0(t),\ldots,X_N(t)):t\ge 0\}$ is defined.
Suppose that this space is equipped with a family of measures $\{\Pm^\nu_x:\nu\in\R_+^{N+1}, x\in\Par^N\}$.
Under $\Pm_x^\nu$, $X_i$ is a Poisson counting process starting at $x_i$ with rate $\nu_i$, and
the elements of $X$ are mutually independent. We set $\Pm^\nu=\Pm^\nu_0$ and
define $T$ as the first exit time from the Weyl chamber $\Par^N$, i.e.,
\[
T=\inf\{t\ge 0: X(t)\not \in \Par^N\}.
\]
The process $X$ killed at time $T$ is a Markov process on $\Par^N$, and
by the Karlin-McGregor formula~\cite{MR0114248} its transition kernel $P^\nu_t$ is given by, for $t>0$,
\begin{equation}
\label{eq:karlinmcgregor}
P^\nu_t(z,z')=
\Pm^\nu_z (X(t)=z',T>t)=\prod_{k=0}^N \left[\nu_k^{z_k'-z_k} e^{-\nu_k t}\right]\det\left\{w_{z'_i-z_j-i+j}(t)\right\}.
\end{equation}
Note that the rates $\nu_i$ do not appear in the determinant,
so that a simple change-of-measure allows us to change jump rates: for $z,z'\in W^N$,
\begin{equation}
\label{eq:RN}
P^\nu_t(z,z') = P^\lambda_t(z,z') \prod_{k=0}^N \left[\left( \frac{\nu_k}{\lambda_k} \right)^{z_k'-z_k}
e^{-(\nu_k-\lambda_k) t}\right].
\end{equation}

For $z,d\in\Par^N$, we define
\begin{eqnarray*}
\Lambda^\nu(z,d)&=&\nu_0^{d_0-z_0}\cdots \nu_N^{d_N-z_N}\det \left\{h^{(jN)}_{z_i-d_j-i+j}\right\}, \\
\Pi^\nu(d,z)&=&\nu_0^{z_0-d_0}\cdots \nu_N^{z_N-d_N}\det \left\{(-1)^{d_i-z_j-i+j}e^{(iN)}_{d_i-z_j-i+j}\right\}.
\end{eqnarray*}
As detailed in Lemma~\ref{lem:Lambda} below, in the terminology of Appendix~A,
$\Lambda^\nu(z,d)$ is a weighted analog of the number of
Gelfand-Tsetlin patterns with specified bottom row and left edge.
The representation of $\phi^\nu_t$ in the following proposition is readily obtained
from the Cauchy-Binet identity (e.g., \cite[Prop.~2.10]{johansson:randommatrices})
\begin{equation}
\label{eq:cauchybinet}
\sum_{z \in \Par^{N}} \det \bigl\{ \xi_i(z_j-j)\bigr\} \det \bigl\{ \psi_j(z_i-i)\bigr\} = \det \left\{
\sum_{z \in \Z} \xi_i(z)\psi_j(z)\right\}.
\end{equation}
For two kernels $A$ and $B$ on $\Par^N$, we define the product kernel $AB$ through
$
AB(z,z')=\sum_{y\in\Par^N} A(z,y) B(y,z').
$

\begin{proposition}
\label{prop:repphi}
We have $\phi^\nu_t=\Pi^\nu P^\nu_t \Lambda^\nu$.
\end{proposition}
\proof
With Cauchy-Binet and the observation
\begin{equation*}
\sum_{r=0}^n (-1)^r e_r^{(iN)}  h^{(jN)}_{n-r}=
\begin{cases}
h^{(ji)}_n    &  \text{ if } j \leq i, \\
(-1)^n e^{(ij)}_n & \text{ if } i\leq j,
\end{cases}
\end{equation*}
we obtain the claim.
\endproof

The kernel $\Pi^\nu$ is the inverse of the kernel $\Lambda^\nu$
as shown in \cite{diekerwarren:determinant2008}.
Proposition~\ref{prop:repphi} is therefore equivalent with
the so-called {\em intertwining} $P^\nu_t\Lambda^\nu=\Lambda^\nu \phi^\nu_t$, which
has the following probabilistic significance. 
There exist processes with values in Gelfand-Tsetlin patterns which give
rise to couplings of the process $D$ and the process $X$ conditioned to stay 
in the Weyl chamber $W^N$.
The kernel $\Lambda^\nu$ describes the conditional
distribution of $D(t)$ given $X(t)$ under this coupling.
Further details can be found in \cite{diekerwarren:determinant2008,warrenwindridge:gelfandtsetlin2008}.

\medskip
The preceding proposition connects departures from queues with non-colliding processes,
and forms the starting point for our investigations.
Another ingredient is the following alternative expression for $\Lambda^\nu$,
which is Proposition 2 of \cite{diekerwarren:determinant2008}. 
It rephrases $\Lambda^\nu$ as a sum over Gelfand-Tsetlin patterns of order $N+1$; 
see Appendix~A for notation and terminology. 

\begin{lemma}
\label{lem:Lambda}
For $z,d\in\Par^N$, we have
\[
\Lambda^\nu( z, d )=\nu_0^{-z_0}\cdots \nu_N^{-z_N}
\sum_{{\mathbf x}\in\mathbf{K}^N: \,\mathrm{sh}(\mathbf{x})=z, \, \mathrm{ledge}({\mathbf x})=d} \nu^{\mathbf{x}}.
\]
\end{lemma}
This lemma shows that  $\Lambda^\nu(z,d)$ vanishes unless $z_N=d_N$, which is not so apparent from the definition.

\section{Transition probabilities for $Q$ and sums of non-crossing probabilities}
\label{sec:identities}
This section presents two identities relating $k^\nu_t(0,0)$ and non-crossing probabilities, where
$k^\nu_t$ is the transition kernel for $Q$: for $q,q'\in\Z^{N}$,
\[
k^\nu_t(q,q')=\pr(Q(t)=q'|Q(0)=q),
\]
suppressing the dependence on $\nu$ on the right-hand side.

Recall that $\Pm^\nu(T>t)$ is the probability that a multidimensional Poisson process
with rate vector $\nu$ stays in the Weyl chamber $\Par^N$ up to time $t$.
We write $\S_{N+1}$ for the symmetric group on $\{0,\ldots,N\}$, i.e.,
all permutations of this set. This group acts on $\R_+^{N+1}$ by permuting the coordinates.
An immediate consequence of the proof of the theorem is that
$k^\nu_t(0,0)$ is symmetric in $\nu_1,\ldots,\nu_N$; this has already been observed
in \cite{oconnell:pathRS2003}.
This property does not extend to $k^\nu_t(q,q')$ for $q,q'\neq 0$.

\begin{theorem}
\label{thm:k00stab}
Suppose that $\nu_1,\ldots,\nu_N$ are distinct. Then we have for $t>0$,
\[
k^\nu_t(0,0)=\sum_{\sigma\in \S_{N+1}:\sigma(N)=0}
\frac {\Pm^{\sigma(\nu)}(T>t)}{\prod_{0\le i<j< N}\left[1-\nu_{\sigma(j)}/\nu_{\sigma(i)}\right]}.
\]
\end{theorem}
\proof
Since $e^{(iN)}_k=0$ for $k<0$ and $i\le N$, it readily follows that 
$\Pi^\nu(0,z)=0$ unless $z=0$, in which case it equals 1. 
Since $h^{(jN)}_k=0$ for $k<0$ and $j\le N$, we similarly find that $\Lambda^\nu(z,d) = 0$ unless $z_N=d_N$,
which can alternatively immediately be deduced from Lemma~\ref{lem:Lambda}. 
On combining Proposition~\ref{prop:repphi} with these facts, we find that for $t>0$,
\begin{eqnarray}
k^\nu_t(0,0)&=&\sum_{\ell\in\Z} \phi^\nu_t(0,(\ell,\ldots,\ell)) = 
\sum_{\ell\in\Z, z\in\Par^N} P^\nu_t(0, z) \Lambda^\nu(z,(\ell,\ldots,\ell))\nonumber \\
&=&\sum_{z\in\Par^N} P_t^\nu(0,z) \Lambda^\nu(z,(z_N,\ldots,z_N)) \\
&=& e^{-\sum_{k=0}^N\nu_k t} \sum_{z\in\Par^N}
\det\left\{w_{z_i-i+j}(t)\right\} \sum_{\mathbf{x}\in\mathbf{K}^N: \mathrm{sh}(\mathbf{x})=z, \mathrm{ledge}(\mathbf{ x})=(z_N,\ldots,z_N)} \nu^{\mathbf x},\nonumber
\end{eqnarray}
where the last equality follows from (\ref{eq:karlinmcgregor}) and Lemma~\ref{lem:Lambda}.
Recall that the relevant notation is introduced in Appendix~A.
Set $\Delta_\nu=\prod_{1\le i<j\le N} (\nu_i-\nu_j)$, a determinant of a Vandermonde matrix.
An empty product should be interpreted as 1.
We next write the sum over $\mathbf x\in\mathbf{K}^N$ in terms of a Schur polynomial:
\begin{eqnarray*}
\sum_{\mathbf{x}\in\mathbf{K}^N: \mathrm{sh}(\mathbf{x})=z, \mathrm{ledge}(\mathbf{ x})=(z_N,\ldots,z_N)} \nu^{\mathbf x} &=& 
\nu_0^{z_N}\cdots \nu_N^{z_N} \sum_{\mathbf{x}\in\mathbf{K}^{N-1}: \mathrm{sh}(\mathbf{x})=(z_0-z_N,\ldots,z_{N-1}-z_N)} 
\nu^{\mathbf x} \\
&=& \nu_0^{z_N}\cdots \nu_N^{z_N} s_{(z_0-z_N,\ldots,z_{N-1}-z_N)}(\nu_1,\ldots,\nu_N) \\
&=& \frac 1{\Delta_\nu} \nu_0^{z_N}\cdots \nu_N^{z_N} \det\left\{\nu_{j}^{z_{i-1}-z_N-i+N}\right\}_{i,j=1,\ldots,N},
\end{eqnarray*}
where the last equality follows from the alternative definition of 
Schur polynomials given in Appendix A.
We thus obtain for $t>0$
\[
k^\nu_t(0,0)=\frac {1}{\Delta_\nu}
e^{-\sum_{k=0}^N\nu_k t}\sum_{z\in \Par^{N}}
\det\left\{w_{z_i-i+j}(t)\right\}
\nu_0^{z_N} \det\left\{\nu_{j}^{z_{i-1}-i+N}\right\}_{i,j=1,\ldots,N}.
\]
The symmetry of $k^\nu_t(0,0)$ with respect to $\nu_1,\ldots,\nu_N$ immediately follows from the symmetry
of the Schur polynomial.
With the Karlin-McGregor formula (\ref{eq:karlinmcgregor}) and the Leibniz formula applied
to the second determinant, we obtain
\[
k^\nu_t(0,0)=\frac {(-1)^N}{\Delta_\nu}\sum_{\sigma\in \S_{N+1}:\sigma(N)=0}\sgn(\sigma)
\nu_{\sigma(0)}^{N-1} \nu_{\sigma(1)}^{N-2}\cdots\nu_{\sigma(N-2)}
\Pm^{\sigma(\nu)}(T>t).
\]
We next note that for $\sigma$ with $\sigma(N)=0$,
\begin{eqnarray*}
\lefteqn{\frac {(-1)^N}{\Delta_\nu} \sgn(\sigma)
\nu_{\sigma(0)}^{N} \nu_{\sigma(1)}^{N-1}\cdots\nu_{\sigma(N-1)}}\\
&=&\frac{\prod_{0\le i<j\le N} \left[\nu_{\sigma(i)}\right]}{(-1)^N \sgn(\sigma) \Delta_\nu} 
=\frac{\prod_{0\le i<j\le N} \left[\nu_{\sigma(i)}\right]}{\prod_{0\le i<j\le N-1} [\nu_{\sigma(i)}-\nu_{\sigma(j)}]} \\
&=&\frac{\prod_{0\le i<j\le N} \left[\nu_{\sigma(i)}\right]\prod_{j=0}^{N-1}[\nu_{\sigma(j)}-\nu_{\sigma(N)}]}{\prod_{0\le i<j\le N} [\nu_{\sigma(i)}-\nu_{\sigma(j)}]} 
=\frac{\prod_{j=0}^{N-1}[\nu_{\sigma(j)}-\nu_{\sigma(N)}]}{\prod_{0\le i<j\le N}\left[1-\nu_{\sigma(j)}/\nu_{\sigma(i)}\right]} \\
&=& \frac{\prod_{j=0}^{N-1}\nu_{\sigma(j)}\prod_{j=0}^{N-1}\left[1-\nu_{\sigma(N)}/\nu_{\sigma(j)}\right]}
{\prod_{0\le i<j\le N}\left[1-\nu_{\sigma(j)}/\nu_{\sigma(i)}\right]}
= \frac{\prod_{j=0}^{N-1}\nu_{\sigma(j)}}
{\prod_{0\le i<j<N}\left[1-\nu_{\sigma(j)}/\nu_{\sigma(i)}\right]},
\end{eqnarray*}
and we have finished the proof.
\endproof

\medskip
For $j=1,\ldots,N$, define the load $\rho_j$ of station $j$ by
$\rho_j=\nu_0/\nu_j$.
We next work under the `stability' assumption $\max_{j=1}^N \rho_j<1$,
and express the difference $k^\nu_t(0,0)-k^\nu_\infty(0,0)$ as a weighted sum of
non-crossing probabilities.

\begin{theorem}
\label{thm:k00rep}
Suppose that $\max_{j=1}^N \rho_j<1$ and that all $\nu_i$ are distinct. Then we have for $t>0$,
\[
k^\nu_t(0,0)=\prod_{j=1}^N\left[1-\rho_j\right]-
\sum_{\sigma\in \S_{N+1}:\sigma(N)\neq 0}
\frac {\prod_{j=1}^N\left[1-\rho_j\right]}{\prod_{0\le i<j\le N}\left[1-\nu_{\sigma(j)}/\nu_{\sigma(i)}\right]}
\Pm^{\sigma(\nu)}(T>t).
\]
\end{theorem}
\proof
For $\lambda\in\R^{N+1}_+$ and $x\in\Par^N$, we define
\[
\omega_\lambda(x)=\lambda_0^{-x_0} \cdots \lambda_N^{-x_N} \det\left\{\lambda_{j}^{x_i-i+j}\right\},
\]
which can be expressed in terms of Schur polynomials (see Appendix~A) as
\[
\omega_\lambda(x) = \prod_{0\le i<j\le N} \left(1-\frac{\lambda_j}{\lambda_i}\right) 
\lambda_0^{-x_0} \cdots \lambda_N^{-x_N}  s_x(\lambda).
\]
For convenience, we write $\omega_\lambda$ for $\omega_\lambda(0)=\prod_{0\le i<j\le N}\left[1-\lambda_{j}/\lambda_{i}\right]$.

Theorem~\ref{thm:k00stab} yields for $t>0$,
\begin{equation}
\label{eq:applyingthm4}
\prod_{j=1}^N\left[1-\rho_j\right]^{-1} k^\nu_t(0,0)=
\sum_{\sigma\in \S_{N+1}:\sigma(N)=0} \omega_{\sigma(\nu)}^{-1} \Pm^{\sigma(\nu)}(T>t).
\end{equation}
Let $\tilde \sigma$ be the permutation for which $\nu_{\tilde \sigma(N)}<\nu_{\tilde \sigma(N-1)}<\ldots<\nu_{\tilde \sigma(0)}$.
The choice of $\tilde \sigma$ entails that the term for $\tilde \sigma$ converges to the constant $\Pm^{\tilde \sigma(\nu)}(T=\infty)$; the idea of the proof is to rewrite
this term.

As outlined in Biane {\em et al.}~\cite[Sec.~5.2]{biane:littlemann2005} for the Brownian case,
we have $\Pm_x^{\tilde \sigma(\nu)}(T=\infty)=\omega_{\tilde \sigma(\nu)}(x)$;
see also \cite[Sec.~4]{oconnell:conditionedRSK2003}.
Moreover, $\omega_{\tilde \sigma(\nu)}(\cdot)$ is harmonic for the killed transition kernel of $X$ under $\Pm^{\tilde \sigma(\nu)}$, meaning that
\[
\sum_{y\in \Par^N} \Pm^{\tilde \sigma(\nu)}_x(X(t)=y,T>t) \omega_{\tilde \sigma(\nu)}(y) =
\omega_{\tilde\sigma(\nu)}(x).
\]
Indeed, using the Karlin-McGregor formula (\ref{eq:karlinmcgregor}), 
this is readily verified with the Cauchy-Binet identity (\ref{eq:cauchybinet}).
We thus obtain
\begin{eqnarray}
\Pm^{\tilde \sigma(\nu)}(T>t)\nonumber
&=&\sum_{x\in \Par^N} \Pm^{\tilde \sigma(\nu)}(X(t)=x,T>t) \left[\omega_{\tilde \sigma(\nu)}(x)+(1-\omega_{\tilde \sigma(\nu)}(x))\right]\nonumber\\
&=& \omega_{\tilde \sigma(\nu)}+ \sum_{x\in \Par^N} \Pm^{\tilde \sigma(\nu)}(X(t)=x,T>t) (1-\omega_{\tilde \sigma(\nu)}(x)).
\label{eq:splitharmonic}
\end{eqnarray}
Next note that the expansion of the determinant $\omega_\nu(x)$ in conjunction with
the change-of-measure formula (\ref{eq:RN}) yield
\begin{eqnarray*}
\lefteqn{\sum_{x\in\Par^N} \Pm^{\nu}(X(t)=x,T>t)
(1-\omega_{\nu}(x))} \\
&=&-\sum_{x\in\Par^N} \Pm^{\nu}(X(t)=x,T>t)
\sum_{\sigma\in \S_{N+1}:\sigma\neq \id} \sgn(\sigma)
\nu_0^{-x_0} \cdots \nu_N^{-x_N} \prod_{i=0}^N \nu_{\sigma(i)}^{x_i-i+\sigma(i)}\\
&=&-\sum_{x\in\Par^N} \sum_{\sigma\in \S_{N+1}:\sigma\neq \id}
\sgn(\sigma) \Pm^{\sigma(\nu)}(X(t)=x,T>t)\prod_{i=0}^N \nu_{\sigma(i)}^{-i+\sigma(i)}\\
&=&-\sum_{\sigma\in \S_{N+1}:\sigma\neq \id}
\sgn(\sigma) \nu_0^0\cdots\nu_N^N \nu_{\sigma(0)}^{-0}\cdots \nu_{\sigma(N)}^{-N}\Pm^{\sigma(\nu)}(T>t)\\
&=&-\sum_{\sigma\in \S_{N+1}:\sigma\neq \id}
\omega_{\nu}\omega_{\sigma(\nu)}^{-1}\Pm^{\sigma(\nu)}(T>t).
\end{eqnarray*}
Replacing $\nu$ by $\tilde \sigma(\nu)$, we see that
the second term in (\ref{eq:splitharmonic}) can be rewritten as
\[
\sum_{x\in\Par^N} \Pm^{\tilde\sigma(\nu)}(X(t)=x,T>t)
(1-\omega_{\tilde\sigma(\nu)}(x)) = -\omega_{\tilde\sigma(\nu)}
\sum_{\sigma\in \S_{N+1}:\sigma\neq \tilde\sigma}
\omega_{\sigma(\nu)}^{-1}\Pm^{\sigma(\nu)}(T>t).
\]
We have thus shown that the right-hand side of (\ref{eq:applyingthm4}) equals
\begin{eqnarray*}
\lefteqn{\sum_{\sigma\in \S_{N+1}:\sigma(N)=0}
\omega_{\sigma(\nu)}^{-1} \Pm^{\sigma(\nu)}(T>t)}\\
&=&\sum_{\sigma\in \S_{N+1}:\sigma(N)=0,\sigma\neq \tilde \sigma}
\omega_{\sigma(\nu)}^{-1}
\Pm^{\sigma(\nu)}(T>t)+1-\sum_{\sigma\in \S_{N+1}:\sigma\neq\tilde \sigma}
\omega_{\sigma(\nu)}^{-1}\Pm^{\sigma(\nu)}(T>t)\\
&=& \sum_{\sigma\in \S_{N+1}:\sigma(N)=0}
\omega_{\sigma(\nu)}^{-1}
\Pm^{\sigma(\nu)}(T>t)+1-\sum_{\sigma\in \S_{N+1}}
\omega_{\sigma(\nu)}^{-1}\Pm^{\sigma(\nu)}(T>t)\\
&=& 1-\sum_{\sigma\in \S_{N+1}:\sigma(N)\neq 0}
\omega_{\sigma(\nu)}^{-1}\Pm^{\sigma(\nu)}(T>t),
\end{eqnarray*}
and the claim follows.
\endproof

\section{Relaxation times}
\label{sec:relaxation}
The representation in Theorem~\ref{thm:k00rep} is particularly suitable
for studying the regime $t\to\infty$. Before giving a result into this direction,
we prove a large-deviation result that we need in our analysis. For terminology, see
Dembo and Zeitouni~\cite{dembozeitouni:largedev1998}.
Recall that the process $X$ is defined on a measurable space equipped with $\Pm^\nu$.
The $\Pm^\nu$-law of $X(t)/t$ satisfies a large-deviation principle
(LDP) in $\R^{N+1}$ on scale $1/t$. The (good) rate function is $I_\nu:\R^{N+1}\to \R$
given by
\[
I_\nu(x)=\begin{cases}
\sum_{k=0}^N \left[x_k \log(x_k/\nu_k) -x_k+\nu_k\right]
&\text{ if } x\in \R_+^{N+1}, \\
\infty
&\text{ otherwise}.
\end{cases}
\]
This can be seen, for instance, as a corollary of Exercise~5.2.12 in \cite{dembozeitouni:largedev1998}.
It turns out that the $\Pm^\nu$-law
of $X(t)/t$ also satisfies an LDP in
\begin{eqnarray*}
\Par_{\R}^N&=&\{x\in\R^{N+1}: x_0\ge \ldots \ge x_N\}
\end{eqnarray*}
on the event that $X$ is killed at the boundary of $\Par^N$.

\begin{lemma}
For any $\nu\in \R_+^{N+1}$, the (defective) measures $\Pm^\nu(X(t)/t\in dx, T>t)$ satisfy
an LDP in $\Par_{\R}^N$ on scale $1/t$. The rate function is $I_\nu$ restricted to $\Par_{\R}^N$.
\end{lemma}
\proof
Since $\Pm^\nu(X(t)/t\in F, T>t)\le \Pm^\nu(X(t)/t\in F)$,
the upper bound follows from the LDP for the $\Pm^\nu$-law of $X(t)/t$, i.e., without the killing mechanism.
For the lower bound, we prove that for any $y$ in the interior of $\Par_{\R}^N$ and any $\delta>0$,
\[
\liminf_{t\to\infty} \frac 1t \log \Pm^\nu \left(X(t)/t\in \prod_{k=0}^N (y_k-\delta,y_k+\delta), T>t\right)\ge -I_\nu(y).
\]
This inequality is trivial if $y\not\in \R_+^{N+1}$, so we assume the contrary.
In view of (\ref{eq:karlinmcgregor}), the probability can be written as
\begin{eqnarray*}
\lefteqn{\int_{\prod_{k=0}^N (y_kt-\delta t,y_kt+\delta t)} \prod_{k=0}^N\left[\left(\frac {\nu_k}{y_k}\right)^{z_k}
e^{-\nu_k t+y_k t}\right]\Pm^y(X(t)\in dz, T>t)}\\
&\ge& \prod_{k=0}^N\left[\left(\frac {\nu_k}{y_k}\right)^{y_kt} e^{-\epsilon |\log(y_k/\nu_k)| t}
e^{-\nu_k t+y_k t}\right]\Pm^y\left(X(t)/t\in \prod_{k=0}^N (y_k-\epsilon,y_k+\epsilon), T>t\right),
\end{eqnarray*}
for any $0<\epsilon<\delta$.
Since the probability on the right-hand side is bounded from below by a positive constant as $t\to \infty$,
we have proven the lower bound after letting $t\to\infty$ and then $\epsilon\to 0$.
\endproof

\medskip
We use the above lemma to investigate the asymptotic behavior of the sum given in Theorem~\ref{thm:k00rep},
as it implies that
\begin{equation}
\label{eq:expdecay}
\lim_{t\to\infty} \frac 1t\log \Pm^{\sigma(\nu)}(T>t)= -\inf_{x\in \Par_{\R}^N} I_{\sigma(\nu)}(x).
\end{equation}
By convexity, the infimum is attained on one of the faces of the polyhedron $\Par_{\R}^N$.
The following result follows from a comparison of the exponential rates of each of the terms
in Theorem~\ref{thm:k00rep}.
We write $\nu_{(i)}$ for the $i$-th smallest rate among $\nu_1,\ldots,\nu_N$, so that
$\nu_{(1)}$ is the smallest rate and $\nu_{(N)}$ the largest. The notation $f(t)\sim g(t)$ is shorthand
for $f(t)/g(t)\to 1$.

\begin{corollary}
\label{cor:asympsurvival}
Suppose that $\max_{j=1}^N \rho_j<1$ and that all $\nu_i$ are distinct.
As $t\to\infty$, we have
\[
\left|k^\nu_t(0,0)-\prod_{j=1}^N \left[1-\rho_j\right]\right|
\sim \frac{\rho_{(1)}}{\prod_{1\le i<j\le N} \left[
1-\nu_{(i)}/\nu_{(j)}\right]} \Pm^{(\nu_{(N)},\ldots,\nu_{(2)},\nu_0, \nu_{(1)})}(T>t).
\]
\end{corollary}
\proof
We may assume that $\nu_0<\nu_1<\ldots<\nu_N$ since $k^\nu_t(0,0)$ is symmetric in $\nu_1,\ldots,\nu_N$.
Theorem~\ref{thm:k00rep} shows that $k_t(0,0)-\prod_{j=1}^N \left[1-\rho_j\right]$ is 
a linear combination of probabilities $\Pm^{\sigma(\nu)}(T>t)$ with $\sigma(N)\neq 0$,
and it is our aim to show that
the term $\Pm^{(\nu_{N},\ldots,\nu_{2},\nu_0, \nu_{1})}(T>t)$ has the slowest exponential decay in $t$.
In view of (\ref{eq:expdecay}), we thus need to prove that $\bar\sigma:=(N,\ldots,2,0,1)$ minimizes
$\inf_{x\in \Par_{\R}^N} I_{\sigma(\nu)}(x)$ over $\sigma\in \S_{N+1}$ with $\sigma(N)\neq 0$.

Let $\sigma$ with $\sigma(N)\neq 0$ be given, and write $j=\sigma^{-1}(0)<N$.
The rate function $I_{\sigma(\nu)}$ is strictly convex on $(0,\infty)^{N+1}$.
Moreover, it achieves its minimum 0 at $x=\sigma(\nu)\not\in W^N_{\R}$;
since $\sigma(\nu)_j< \sigma(\nu)_N$, the latter point is
separated from $W^N_{\R}$ by the hyperplane $\{x\in\R^{N+1}: x_j=x_N\}$.
The infimum of $I_{\sigma(\nu)}$ over $W^N_{\R}$ is thus at least as large as the infimum over this
hyperplane, with equality if and only if the minimizing argument lies in $W^N_{\R}$.
It is not hard to see that the infimum over the hyperplane is $\nu_0 +\nu_{\sigma(\nu)_N} - 2\sqrt{\nu_0 \nu_{\sigma(\nu)_N}}$,
and that it is achieved at $\nu^*(\sigma)$ given by
\[
\nu^*(\sigma)_i=
\begin{cases}
\sigma(\nu)_i & \text{ if } i\not \in \{j, N\}, \\
\sqrt{\nu_0\sigma(\nu)_N} & \text{ otherwise}.
\end{cases}
\]
We conclude that, if $j=\sigma^{-1}(0)<N$,
\begin{eqnarray*}
\inf_{x\in W^N_\R} I_{\sigma(\nu)}(x) &\ge&
\inf_{x\in (0,\infty)^{N+1}: x_j=x_N} I_{\sigma(\nu)}(x) \\&=&
\nu_0+\nu_{\sigma(\nu)_N}-2\sqrt{\nu_0\nu_{\sigma(\nu)_N}}\\ 
&\ge& \nu_0+\nu_1-2\sqrt{\nu_0\nu_1},
\end{eqnarray*}
where the last equality uses that $x\mapsto \nu_0 + x - 2\sqrt{x\nu_0}$ is increasing in $x\ge \nu_0$.

We next argue that this lower bound is attained only for $\sigma=\bar\sigma$.
By strict convexity of $I_{\sigma(\nu)}$,
the first inequality in the preceding display is strict unless $\nu^*(\sigma)\in W^N_\R$, showing
that $\sigma=\bar \sigma$ is the only candidate to attain the lower bound.
In fact, since both inequalities become equalities if $\sigma=\bar\sigma$,
$\bar \sigma$ corresponds to
the unique asymptotically dominant term in Theorem~\ref{thm:k00rep}.

A direct calculation shows that for $\sigma=\bar\sigma$,
\[
\frac {\prod_{j=1}^N\left[1-\rho_j\right]}{\prod_{0\le i<j\le N}\left[1-\nu_{\sigma(j)}/\nu_{\sigma(i)}\right]}
= \frac{-\rho_1}{\prod_{1\le i<j\le N} \left[
1-\nu_{i}/\nu_{j}\right]},
\]
and this determines the prefactor.
\endproof

\medskip
The relaxation time is defined as the reciprocal of the exponential rate at which
$|k^\nu_t(0,0)-k^\nu_\infty(0,0)|$ converges to zero.
The proof of Corollary~\ref{cor:asympsurvival} shows in particular that, 
if the $\nu_i$ are all distinct and the system is stable, the relaxation time is
\begin{equation}
\label{eq:relaxationtime}
\frac1{\nu_0+\nu_{(1)}-2\sqrt{\nu_0\nu_{(1)}}}=
\frac{1}{\nu_{(1)}(1-\sqrt{\rho_{(1)}})^2}.
\end{equation}
This is the same as the relaxation time of an M/M/1 queue with arrival rate $\nu_0$ and service rate $\nu_{(1)}$.
Therefore, the relaxation time is determined by a single {\em bottleneck} station.
These observations prove a conjecture of Blanc~\cite{blanc:relaxation1985}.

Corollary~\ref{cor:asympsurvival} is also of interest in order to obtain {\em exact} asymptotics (as $t\to\infty$)
for the absolute difference between $k^\nu_t(0,0)$ and $k^\nu_\infty(0,0)$,
in contrast with the rough (i.e., logarithmic) asymptotics obtained here.
In fact, exact tail asymptotics for non-crossing probabilities are a current research topic.
For recent work on such asymptotics in a Brownian setting, refer to
Pucha{\l}a and Rolski~\cite{puchalarolski:asymptotics2007}.
Thus, a Poisson-process analog of the results in \cite{puchalarolski:asymptotics2007}
would immediately yield new results for queues in series.

\section{Generalizations}
\label{sec:discussion}
This section discusses analogs of the identity presented in Theorem~\ref{thm:k00stab} for $k_t^\nu(q,q')$ with
$q,q'\neq 0$, and how it changes if the service rates are not necessarily distinct.

Define $\pi:W^N \rightarrow \Z^N_+$ via $\pi(x)=(x_0-x_1, x_1-x_2, \ldots x_{N-1}-x_{N})$ and recall that the queue length vector and departures vector are related by $Q(t)= Q(0)+ \pi(D(t))$. For $\ell\in\Z$, we also define the `inverse'
$\invpi_\ell: \Z^N_+\rightarrow W^N$ of $\pi$ through $\invpi_\ell(q) =
(\ell+q_1+\ldots+q_N, \ell+q_2+\ldots+q_{N},\ldots,\ell+q_N,\ell)$, and we abbreviate $\invpi_0$
by $\invpi$.

The transition probabilities for $Q$ can be expressed in terms of those of $D$ through
$k^\nu_t(q, q')= \sum_{\ell\in\Z} \phi^\nu_t( \invpi(q), \invpi_\ell(q'))$, for $q,q' \in {\mathbb Z}^N_+$.
Now if we define for $z\in W^N$, $q\in\Z_+^{N}$,
\begin{eqnarray}
\Lambda^\nu_Q(z,q)&=&\sum_{\ell\in\Z}\Lambda^\nu(z,\invpi_\ell(q))=
\Lambda^\nu(z,\invpi_{z_N}(q))\nonumber\\
&=&\nu_0^{-z_0}\cdots\nu_N^{-z_N}
\sum_{{\mathbf x}\in\mathbf{K}^N: \,\mathrm{sh}(\mathbf{x})=z, \,
\mathrm{ledge}({\mathbf x})=\invpi_{z_N}(q)} \nu^{\mathbf{x}}\label{eq:defLambdaQ} \\
&=&\prod_{k=0}^{N-1} \left[\nu_k^{\invpi(q)_k-z_k}\right]
\det\left\{\nu_j^{z_N} h^{(jN)}_{z_i-z_N-\invpi(q)_j
-i+j}\right\}_{i,j=0,\ldots,N-1},\nonumber
\end{eqnarray}
then it follows from $P^\nu_t\Lambda^\nu=\Lambda^\nu \phi^\nu_t$ that
$P^\nu_t$ and $k^\nu_t$ are intertwined: 
\[
P_t^\nu\Lambda_Q^\nu(z,q)=\sum_{\ell\in\Z}P_t^\nu
\Lambda^\nu(z,\invpi_\ell(q))=\sum_{\ell\in\Z}
\Lambda^\nu\phi_t^\nu (z,\invpi_\ell(q))=
\Lambda_Q^\nu k_t^\nu(z,q).
\]
The kernel $\Lambda_Q^\nu$ is not invertible, however it does admit a (non-unique) left inverse $\Pi_Q^\nu$ defined by  $\Pi_Q^\nu(q,z)=\Pi^\nu(\invpi(q),z)$: for $q,q'\in\Z^N_+$,
\[
\Pi_Q^\nu \Lambda_Q^\nu(q,q') = \sum_{\ell\in\Z} \sum_{z\in W^N} \Pi^\nu(\invpi(q),z) 
\Lambda^\nu(z,\invpi_\ell(q')) = \sum_{\ell\in\Z} 1_{\{\invpi(q)=\invpi_\ell(q')\}} = 1_{\{q=q'\}}.
\]
Consequently we obtain the representation of the queue length transition probabilities, which also follows
directly from Proposition~\ref{prop:repphi}:
\begin{proposition}
\label{prop:kernelQ}
For any $t>0$, we have
\[
k^\nu_t= \Pi_Q^\nu P^\nu_t \Lambda_Q^\nu.
\]
\end{proposition}

The above proposition may equivalently be written as
\begin{equation}
\label{eq:repexpectedval}
k^\nu_t(q,q^\prime)= \sum_{z\in W^N} \Pi_Q^\nu (q,z)
\meanm^\nu_{z} \left[ \Lambda_Q^\nu(X(t),{q^\prime}), T>t \right].
\end{equation}
The function $z\mapsto \Pi_Q^\nu(q,z)$ is supported on a finite set for fixed $q$, as is readily verified
from its definition. Therefore,
the sum over $z\in W^N$ in (\ref{eq:repexpectedval}) is in fact a finite sum.

Of special interest is the case $q^\prime=0$. 
From (\ref{eq:defLambdaQ}) it then follows that $\Lambda_Q^\nu(z,q^\prime)$ can be expressed in terms of a Schur polynomial for nonnegative $z\in W^N$: 
\begin{eqnarray*}
\Lambda^\nu_Q(z,0)&=& 
\nu_0^{-z_0}\cdots\nu_N^{-z_N}
\sum_{{\mathbf x}\in\mathbf{K}^N: \,\mathrm{sh}(\mathbf{x})=z, \,
\mathrm{ledge}({\mathbf x})=(z_N,\ldots,z_N)} \nu^{\mathbf{x}} \\&=& 
\nu_0^{z_N-z_0}\cdots\nu_{N-1}^{z_N-z_{N-1}} s_{(z_0-z_N,\ldots,z_{N-1}-z_N)}(\nu_1,\ldots,\nu_N).
\end{eqnarray*}
It is only if the rates $\nu_i$ are distinct that this can be expressed as a sum of geometric terms and a result of the form of Theorem~\ref{thm:k00stab} can be obtained from (\ref{eq:repexpectedval}) by suitable changes of measure, cf.~the determinantal representation of Schur polynomials in Appendix~A.
At the other extreme, when all the $\nu_i$ are equal to $\nu$, we have
$\Lambda_Q^\nu(z,0)=\prod_{0\le i<j\le N-1} (z_i-z_j-i+j)$ and (\ref{eq:repexpectedval}) becomes
\[
k^\nu_t(q,0)= \sum_{z\in W^N} \Pi_Q^\nu (q,z)
\meanm^\nu_{z} \left[ \prod_{0\le i<j\le N-1} (X_i(t)-X_j(t)-i+j), T>t \right].
\]

If $q^\prime \neq 0$ and the rates $\nu_i$ are distinct, then $\Lambda_Q^\nu(z,q^\prime)$ has a piecewise
sum-of-geometric form as long as $z\ge \invpi(q^\prime)$, otherwise it is zero.
We expect that Theorem~\ref{thm:k00stab} can be generalized into this direction, but that it involves the sum of
probabilities of the form $\{X(t) \ge \invpi(q^\prime), T>t\}$.

We close this section by relating Proposition~\ref{prop:kernelQ} to the literature.
In the single-station case $N=1$, we have $\Lambda_Q^\nu(z,q)= (\nu_1/\nu_0)^{z_0-z_1-q} 1_{\{z_0\ge z_1+q\}}$
and $\Pi_Q^\nu(q,z)=1_{\{z=(q,0)\}} - \nu_1 1_{\{z=(q-1,0)\}}$ for $q\ge 1$ while $\Pi_Q^\nu(0,z)=1_{\{z=0\}}$.
This leads to
\[
k^\nu_t(q,q') = e^{-(\nu_0+\nu_1)t}
\mathop{\sum_{z\in \Par^1}}_{z_0\ge z_1+q'} \nu_0^{z_1+q'-q} \nu_1^{z_0-q'}
\left|\begin{array}{cc}w_{z_0-q}(t)-\nu_1 w_{z_0-q+1}(t) & w_{z_0+1}(t)\\w_{z_1-q-1}(t)-\nu_1 w_{z_1-q}(t)& w_{z_1}(t)\end{array}\right|.
\]
We know from the sum representation of the modified Bessel function (e.g., \cite[9.6.10]{abramowitzstegun}) that
\[
\sum_{k\in\Z} \nu_0^k \nu_1^k w_{k+d}(t) w_k(t)= (\nu_0\nu_1)^{-d/2} I_d(2\sqrt{\nu_0\nu_1}t)
\]
for $d\in\Z$. Expanding the determinant, we find that $k^\nu_t(q,q')$ equals
\[
\rho^{(q'-q)/2} I_{q'-q} + \rho^{(q'-q-1)/2} I_{q+q'+1} + (1-\rho)\rho^{q'} \sum_{\ell\ge q+q'+2} \rho^{-\ell/2} I_\ell,
\]
where $\rho:=\nu_0/\nu_1$ and all arguments of the Bessel functions are $2\sqrt{\nu_0\nu_1} t$.
This is the well-known expression for the time-dependent M/M/1 queue.
A similar program can be followed for $N>1$,
and $k^\nu_t(q,q')$ can then be interpreted as weighted sum
of so-called lattice Bessel functions \cite{baccellimassey:transient1988,massey:exit1987,bohmmohanty:karlinmcgregor1997}.
To our knowledge, for $N>2$ this representation has not been recorded in the literature.
Note that it is numerically inefficient to use the resulting representation of $k^\nu_t$,
see \cite{bohmjainmohanty:combinatorial1993}.

\section*{Acknowledgments}
A.B.~Dieker was financially supported in part by a postdoctoral fellowship from the IBM T.J.~Watson Research Center, Yorktown Heights NY, USA.

\appendix
\section{Symmetric functions and Gelfand-Tsetlin patterns}
This appendix defines the coefficients $e_k^{(ij)}$ and $h_k^{(ij)}$ used in Section~\ref{sec:departures}.
We also introduce Gelfand-Tsetlin patterns and Schur polynomials, which play important roles in the proof
of our main results.
More details can be found in, e.g., Chapter 7 of Stanley~\cite{stanley:ec2}.

The $r$th complete homogeneous symmetric polynomials in the indeterminates $\alpha_0,\ldots
\alpha_N$ is given by
\[
h_r(\alpha)=\sum_{k_0\geq 0, \ldots, k_N \geq 0: k_0+k_1+\cdots +k_N=r}
\alpha_0^{k_0}\alpha_1^{k_1} \cdots \alpha_N^{k_N}.
\]
By convention $h_0=1$ and $h_r=0$ for $r<0$.  Now for $0\leq i <
j \leq N$, let $h^{(ij)}_r(\alpha) = h_r(\alpha^{(ij)})$ where
$\alpha^{(ij)}$ is the $(N+1)$-vector
$(0,\ldots,0,\alpha_{i+1},\alpha_{i+2}, \ldots ,\alpha_j,0, \ldots
0)$ obtained from $\alpha$ by setting the first $i+1$ weights and the
last $N-j$  weights  equal to $0$. Equivalently it is the $r$th
complete homogeneous symmetric polynomial in the indeterminates
$\alpha_{i+1}, \ldots, \alpha_j$.
We set $h^{(jj)}_r(\alpha)= {\mathbf 1}(r=0)$.

We write $e_r$ for the $r$th elementary symmetric polynomial defined as 
\[
e_r(\alpha)= \sum _{0\le k_1<k_2< \cdots <k_r\le N} \alpha_{k_1} \cdots \alpha_{k_r}.
\]
In analogy with the complete homogeneous symmetric polynomials, we
use the conventions $e^{(jj)}_r(\alpha)={\mathbf 1}(r=0)$ and $e_0(\alpha)=1$. 
We also set $e^{(ij)}_r(\alpha)= e_r(\alpha^{(ij)})$, so that in particular $e_r^{(ij)}(\alpha)=0$ for $r<0$.

Let $\mathbf{x}$ be an array of real-valued variables
$\mathbf{x}=(x^0,\ldots,x^N)$ with $x^k=(x_0^k,x_1^k,\ldots,x_k^k)\in {\mathbf Z}^{k+1}$,
such that the coordinates satisfy the inequalities
\[
x^k_k\le x_{k-1}^{k-1}\leq x^k_{k-1}\leq x^{k-1}_{k-2}
\le \ldots\leq x_1^k\leq x_0^{k-1}\leq x_0^k
\]
for $k=1,\ldots,N$.
We write $\mathbf{K}^N$ for the set of all $\mathbf x$ satisfying the above constraint,
and say that any $\mathbf x \in \mathbf{K}^N$ is a {\em Gelfand-Tsetlin (GT) pattern} of order $N+1$.
For $\mathbf{x}\in\mathbf{K}^N$, we set
$\mathrm{sh}(\mathbf{x})=(x^N_0,x^N_1, \ldots, x_N^N)$ and
$\mathrm{ledge}(\mathbf{ x})= (x^0_0,\ldots,x^N_N)$.
If $\mathrm{ledge}(\mathbf{ x})\ge 0$, i.e., all $x_i^j$ are nonnegative, 
integer-valued GT patterns of order $N+1$ are in one-to-one correspondence with so-called
semistandard Young tableau with $N+1$ rows and entries not exceeding $N+1$. 
For a vector $\alpha=(\alpha_0,\alpha_1, \ldots, \alpha_N)$ of weights,
we define the weight $\alpha^{\mathbf{x}}$ of a GT pattern $\mathbf{x}$ by
\[
\alpha^{\mathbf{x}}= \alpha_0^{x^0_0} \prod_{k=1}^N \alpha_k^{ \sum
x^k_i- \sum x^{k-1}_i}.
\]

The Schur polynomial in the indeterminates $\alpha_0,\ldots
\alpha_N$ corresponding to a nonnegative $z\in\Par^N$ is given by
\[
s_z(\alpha)= \sum_{\mathbf x \in \mathbf{K}^N: \mathrm{sh}(\mathbf{x})=z}
\alpha^{\mathbf x}.
\]
Schur polynomials are symmetric in the $\alpha_i$, which is readily seen from their alternative determinantal definition:
\[
s_z(\alpha)=\frac 1{\prod_{0\le i<j\le N} [\alpha_i-\alpha_j]} \det\left\{ \alpha_j^{z_i-i+N}\right\}.
\]
This representation requires that the $\alpha_i$ be distinct,
but there is no singularity; 
there are factors of $(\alpha_i-\alpha_j)$ implicit in the determinant which cancel with factors in the
denominator.

{\small
\bibliography{../../../bibdb}
\bibliographystyle{amsplain}
}

\end{document}